\input epsf.sty
\input amstex
\documentstyle{amsppt} 
\magnification=1200
\hcorrection{.25in}
\advance\voffset.25in
\advance\vsize-.5in

\topmatter \title Sierksma's Dutch Cheese Problem \endtitle \author K. S. Sarkaria\endauthor
\address \eightpoint K. S. Sarkaria, 
 Department of Mathematics, 
 Panjab University, 
 Chandigarh 160014, INDIA. \endaddress
\thanks Part of this work was done at the Mathematical Sciences
Research Institute, Berkeley.
Research at MSRI is supported in part by NSF grant DMS-9022140.\endthanks

\abstract Consider partitions of a cardinality $(q-1)(d+1)+1$ generic subset of euclidean $d$-space, into $q$ parts whose convex hulls have a nonempty intersection. We show that if these partitions are counted with appropriate signs $\pm 1$ then the answer is always $((q-1)!)^d$. Also some other related results are given.\endabstract 

\endtopmatter
\NoBlackBoxes

\document 

\heading1. Introduction. \endheading The object of this note is to prove the following, thus establishing a conjecture of Sierksma [10], 1979.

\medpagebreak

\bf Theorem 1. \it Let $S$ be any cardinality $(q-1)(d+1)+1$ subset of a real affine $d$-dimensional space ${\Bbb A}^d$. Then there exist at least $((q-1)!)^d$ partitions of $S$ into $q$ disjoint subsets, $S = \sigma_1 \cup \dots \cup \sigma_q$, with \rm $\text{conv}(\sigma_1)\cap \dots \cap \text{ conv}(\sigma_q)$ \it  nonempty.

\medpagebreak

\rm We remark that even the existence of \it one \rm such partition is not obvious : it was conjectured by Birch [2] in 1958 and confirmed by Tverberg [13] in 1966 only by means of a fairly involved argument (but see also remark (e) of \S 5). Henceforth we will refer to partitions of the above kind as \it Tverberg partitions \rm of $S$.

The proof of the above theorem is given in \S 4 and depends on Lemma 1
of \S 2 which verifies that the Euler number of a certain complex
vector bundle $\frak L^{\perp}$ on complex projective
$(q-1)(d+1)$-space is $((q-1)!)^{d+1}$, and on Lemmas 2-4 of \S 3
which serve to relate the Tverberg partitions of $S$ with the zeros of
a section $s$ of $\frak L ^{\perp}$. (For basic facts regarding
characteristic classes see Steenrod [11], Hirzebruch [3], and
Milnor-Stasheff [5].) We in fact arrive at a precise \it index formula
\rm which shows that if one counts the Tverberg partitions of a
generic $S$ with appropriate signs $\pm$ then the answer is always
$((q-1)!)^d$.  We note in \S 5 that our method also establishes the
so-called \it \lq\lq continuous\rq\rq \ Tverberg theorem for all $q$.

\smallpagebreak

The bound given by the above theorem is the best possible \rm as can
be seen by using Sierksma's configuration $S_0$ : take $q-1$
coincident points at each of $d+1$ affinely independent positions, and
let the last point be at the barycenter (the case $d=2$, $q=3$ is
shown below).

$$\epsffile{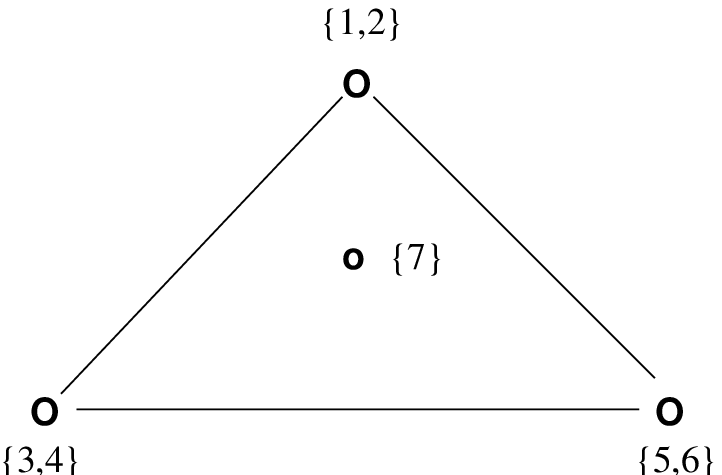}$$

To obtain a Tverberg partition of $S_0$, totally order each set of
coincident points, and for each $1 \leq i \leq q-1$, let $\sigma_i$
consist of the $d+1$ points which occur at the $i$th places in these
orderings, and let $\sigma_q =$ \{barycenter\}. The number of these
partitions is $((q-1)!)^d$ and it is clear that there is no other
Tverberg partition of $S_0$.

\heading 2. A Chern number \endheading 

We denote by $\Bbb Z /q$ the group of the $q$th roots of unity $$\{ 1,
\omega = \text{exp} (2 \pi i/q), \omega ^2 , \dots , \omega ^{q-1}
\},$$ and consider the \it regular representation $ \Bbb V$ \rm of this
group, i.e. the vector space of all $q$-tuples of complex numbers
equipped with the action of $\omega$ given by
$$ (z_1, z_2, \dots , z_q) \longmapsto (z_2, \dots , z_q, z_1).$$

Now $\Bbb V$ is a direct sum $\oplus_{0 \leq k < q} {\Bbb V}_k$ of q
inequivalent one dimensional representations, with the $k$th being
given by

$$ {\Bbb V}_k \ = \ \{ (z_1, z_2, \dots , z_q) : z_{i+1} = {\omega}^k
. z_i \forall i \}.$$

\flushpar We note that the diagonal ${\Bbb V}^0 = \Bbb V_0$ is given by $z_1 = \dots = z_q$, while its complement ${\Bbb V}^{\perp} = \oplus _{1 \leq k <q} {\Bbb V}_k$ is given by $z_1 + z_2 + \dots + z_q = 0$.

We note next that ${\Bbb Z}/q$ is a subgroup of the group $S^1$ of
complex numbers of absolute value 1, and that its action on $\Bbb V
_k$ extends in a natural way to an action of $S^1$ (or even of $\Bbb C
^ \times$) on $\Bbb V _k$ :

$$g \bullet (z_1, z_2, \dots , z_q) \ = \ (g^k z_1, g^k z_2, \dots ,
g^k z_q)\ \ \forall g \in S^1.$$

Using this we'll consider each $\Bbb V _k$, and so also their direct
sum $\Bbb V$, as a representation of $S^1$ (or even of $\Bbb C
^\times$), and the $(d+1)$-fold direct sum $\Bbb V \oplus \dots \oplus
\Bbb V$ will be denoted $\Bbb L$. Likewise $\Bbb L ^0 = \Bbb V ^0
\oplus \dots \oplus \Bbb V ^0$ and $\Bbb L ^{\perp} = \Bbb V ^\perp
\oplus \dots \oplus \Bbb V ^\perp$.

We now consider the \it Hopf bundle \rm $\xi$, i.e. the principal
$S^1$-bundle over $\Bbb C P^N$, $N = (q-1)(d+1)$, whose total space is
the sphere

$$ S^{2N+1} = \{ (z_1, \dots ,z_{N+1}) \in {\Bbb C}^{N+1} : |z_1| +
\dots + |z_{N+1}| = 1 \}$$

\flushpar of ${\Bbb C}^{N+1}$, with the action of $S^1$ being given by multiplication.  We denote by $\frak L$ the quotient of $S^{2N+1} \times \Bbb L$ by the diagonal $S^1$-action, i.e. the total space of the complex vector bundle $\frak L \rightarrow \Bbb C P^N$ associated to $\xi$ by the $S^1$-representation $\Bbb L$.

\smallpagebreak

\bf Lemma 1. \it If $c_N(\frak L)$ denotes the $N$th Chern class of $\frak L$, then
$$ < c_N(\frak L), \Bbb C P^N >\ = \ ((q-1)!)^{d+1}.$$

Proof. \rm We will use the fact that the first Chern class $c_1(\xi )
\in H^2(\Bbb C P^N ; \Bbb Z)$ generates the cohomology ring $H^*(\Bbb
C P^N ; \Bbb Z)$, i.e. that $H^* (\Bbb C P^N ; \Bbb Z)$ consists of
all integral multiples of its powers $(c_1(\xi ))^j$.

Let $\frak V _k$ be the complex line bundle associated to $\xi$ by the
irreducible representation $\Bbb V_k , 0 \leq i < q$.  We note that
$\frak L$ is a direct sum of $d+1$ copies of each of these $q$ line
bundles (so it has complex fibre dimension $N+(d+1)$).  Thus
$c_N(\frak L) \in H^{2N}(\Bbb C P^N ; \Bbb Z)$ is the $(d+1)$th power
of the cohomology class

$$c_1(\frak V _1)\cdot c_1(\frak V _2)\cdot \dots \cdot c_1(\frak
V_{q-1}) \in H^{2(q-1)}(\Bbb C P^N ; \Bbb Z).$$

Thus what we need to check is that this is $(q-1)!$ times the
generator $(c_1(\xi ))^{q-1}$ of $H^{2(q-1)}(\Bbb C P^N ; \Bbb Z )
\cong \Bbb Z$.  This follows because $c_1 (\xi ) = c_1 (\frak V _1)$
and $c_1(\frak{ V} _k \cong \otimes ^k \frak V _1 ) = k.c_1(\xi)$. \it
q.e.d.\rm

We note that $c_N(\xi)$ is also the \it Euler class \rm of the vector
bundle $\frak L ^\perp \rightarrow \Bbb C P^N$ associated to $\xi$ by
the sub representation $\Bbb L^\perp$, and what we have calculated is
the \it Euler number \rm of this oriented vector bundle.

\heading 3. A deleted join
\endheading  

We will denote by $K$ the simplicial complex consisting of all subsets
of the standard basis $\{ e_1, \dots , e_{N+1}\}$ of $\Bbb C^{N+1}$.
The sphere $S^{2N+1}$ is the join of the circles $ge_\alpha$, $g\in
S^1$, i.e. it consists of all points $P = \sum_\alpha t_\alpha
(P)g_\alpha (P)e_\alpha$. Here (and in all similar sums) $\alpha$ runs
over some subset $[P]$ of cardinality $|P|$ of $\{ 1,2, \dots , N+1
\}$, $g_\alpha (P) \in S^1$, and the $t_\alpha(P)$'s are positive
reals with sum 1. (Frequently we'll drop the $(P)$'s and simply write
$ \sum_\alpha t_\alpha g_\alpha e_\alpha$ etc.) Collecting together
terms having equal $g_\alpha$'s we will also sometimes write $P =
\sum_{i=1}^{i=r(P)}t_i(P)g_i(P)x_i(P)$ where $1\leq r(P) \leq |P| \leq
N+1$ denotes the number of \it distinct \rm $g_\alpha(P)$'s, $t_i(P) =
\sum_{g_\alpha = g_i}t_\alpha (P)$ 's are positive numbers having sum
1, and $x_i(P) = \frac{\sum _{g_\alpha = g_i}t_\alpha
(P)e_\alpha}{\sum_{g_\alpha = g_i}t_\alpha (P)}$'s are points
belonging to \it pairwise disjoint \rm faces of the geometrical
simplex $|K| = \text{conv} \{e_1, \dots , e_{N+1} \}$.

Recall now that the simplices $(\sigma _1, \sigma _2, \dots , \sigma
_q)$ of the join $K\cdot K \cdot \dots \cdot K$ of $q$ disjoint copies
of $K$ are obtained by taking unions of $q$ simplices one from each
copy.  We equip it with the $\Bbb Z/q$-action $(\sigma _1, \sigma _2 ,
\dots , \sigma _q) \mapsto (\sigma _2, \dots , \sigma _q, \sigma _1)$.
The $q$-fold \it deleted join \rm $K\ast K \ast \dots \ast K$ is the
free $\Bbb Z/q$-subcomplex of $K\cdot K\cdot \dots \cdot K$ consisting
of all simplices $(\sigma _1, \sigma _2, \dots , \sigma _q)$ for which
the $\sigma$'s are disjoint in $K$.

Taking the $j$th copy, $1 \leq j \leq q$, of $K$ to be that having
vertices $\omega ^{j-1}e_\alpha$, we will identify the geometrical
realization of $K\ast \dots \ast K$ with the free $\Bbb Z/q$-subspace
of $S^{2N+1}$ consisting of all points $P =
\sum_{i=1}^{i=r(P)}t_i(P)g_i(P)x_i(P)$ with $g_i(P) \in \Bbb Z/q
\subset S^1$ (so $r(P) \leq q$ here).  We'll denote by $K\# K \#
\cdots \# K$ the $q$-fold \it deleted product \rm of $K$, i.e. the
subspace of $| K\ast \dots \ast K|$ consisting of all points $P$ of
the form $\sum_{i=1}^{q}\frac{1}{q}g_i(P)x_i(P)$, $g_i(P) \in \Bbb
Z/q$.

We will consider $\Bbb L = \Bbb V \oplus \dots \oplus \Bbb V$ as \it
the vector space $\Bbb C ^{d+1} \otimes \Bbb V$ of all $(d+1) \times
q$ matrices \rm over $\Bbb C$, and $\Bbb A ^d$ (\it resp. $\Bbb
R^{d+1}$, resp. $\Bbb C^{d+1}$\rm ) \rm as the affine (\it resp. \rm
linear, \it resp. \rm complex linear) span in $\Bbb L$ of the $d+1$
matrices whose sole nonzero entry is 1 and lies in the first column.
A matrix will often be denoted by the sequence of its column vectors.

The given points $S = \{ s_1, \dots ,s_{N+1} \}$ of $\Bbb A^d \subset
\Bbb L$ determine the affine linear map $|K| @>s>> \Bbb A^d$,
$s(e_\alpha) = s_\alpha$, which we'll extend to the continuous
$S^1$-map $s: S^{2N+1} \rightarrow \Bbb L$ defined by $s(P) =
\sum_\alpha t_\alpha (P)(g_\alpha (P) \bullet s_\alpha )$. The
decomposition $\Bbb L = \Bbb L ^0 \oplus \Bbb L^{\perp}$ then
determines an $S^1$-map $s:S^{2N+1} \rightarrow \Bbb L^{\perp}$ whose
zeros are the principal items of interest for us.

\medpagebreak

\bf Lemma 2. \it Any zero $P$, $s(P)=0$, of the $S^1$-map $s :S^{2N+1} \rightarrow \Bbb L ^\perp$ must either lie on an $S^1$-orbit passing through the $q$-fold deleted product $K\# \cdots \# K$, or else be such that $r(P) > q$ with the projections of the points $\{ g_\alpha (P)\bullet s_\alpha : \alpha \in [P] \}$ \rm not \it in general position on $\Bbb L^\perp$. \rm (Here  as usual \lq\lq general position\rq\rq \  means that any $i \leq $ dimension of the vectors are linearly independent.) 

\medpagebreak

\it Proof. \rm Using the definitions of the matrix space $\Bbb L$ and the $S^1$-map $s : S^{2N+1} \rightarrow \Bbb L$ we have
$$
s(\sum_\alpha t_\alpha g_\alpha e_\alpha) =
\frac{1}{q}\sum_{l=0}^{q}\{\sum_\alpha t_\alpha g_\alpha ^{\
l}s_\alpha (1,\omega ^l, \dots ,\omega ^{l(q-1)})\}.$$

The orbit of $P = \sum_{\alpha} t_{\alpha}g_\alpha e_\alpha \in
S^{2N+1}$ is mapped by $s$ to a single point iff the right side of the
above equation is in $\Bbb L^0$, i.e. iff only the first of its
summands is nonzero, i.e. iff

$$ \sum _{\alpha = 1}^{N+1} t_\alpha {g_\alpha}^ls_\alpha = 0,\ \ \ 1
\leq l \leq q-1.$$

Alternatively, using $P = \sum_{i=1}^r t_i g_i x_i $ (recall that here
the $g_i$'s are distinct), these equations read

$$\sum_{i=1}^r t_i{g_i}^ly_i = 0, \ \ 1 \leq l \leq q-1,$$

\flushpar where $y_i = s(x_i)$.

\smallpagebreak

CASE $r(P) < q$.  This cannot happen because now the equations
$\sum_{i=1}^r z_i{g_i}^l = 0$, $1 \leq l \leq q-1$, have only the
trivial solution, whereas $t_iy_i \in \Bbb R^{d+1} \backslash \{ 0
\}$.

\smallpagebreak

CASE $r(P) = q$. We will denote the Vandermonde determinant
$|{g_i}^l|$, $1 \leq l \leq q-1$, $1 \leq i \leq q-1$, by $\Delta
(g_1, \dots ,g_{q-1})$ or just $\Delta$.  It is related to its complex
conjugate by

$$ \overline{\Delta (g_1, \dots ,g_{q-1})} = (-1)^{\binom{q-1}{2}}(g_1
\cdots g_{q-1})^{-q} \Delta (g_1, \dots ,g_{q-1}).$$

\flushpar This follows from   

$$\Delta (g_1, \dots ,g_{q-1}) = g_1 \cdots g_{q-1} \prod_{i>j} (g_i -
g_j)$$

\flushpar and the fact that the complex conjugate of any $g\in S^1$ is the same as its inverse.

By Cramer's rule any solution of the system of equations $\sum_{i=1}^q
z_i {g_i}^l = 0$, $1 \leq l \leq q-1$, is of the form $z_m = - \frac
{\Delta_{m,q}}{\Delta}z_q$, $1 \leq m \leq q$, with $z_q$
arbitrary. Here $\Delta_{m,q}$ denotes the determinant obtained from
$\Delta$ by replacing $g_m$ by $g_q$. \it These determinantal ratios
are all real only if $\{ g_1, \dots ,g_q \}$ is a coset of $\Bbb Z/q$
in $S^1$. \rm This is so because by above the complex conjugate of
$\frac {\Delta_{m,q}}{\Delta}$ is $(\frac{g_m}{g_q})^q\frac
{\Delta_{m,q}}{\Delta}$ and so for reality we must have
$(\frac{g_m}{g_q})^q = 1$.

Since $t_iy_i \in \Bbb R^{d+1}$ it follows that $\{ g_1, \dots ,g_q \}
= g^{-1}\cdot \Bbb Z/q$ for some $S^1$.  In this case any solution of
our system of linear equations satisfies $z_1 = \dots = z_q$, so the
$t_iy_i$'s must be equal to each other, which implies that the points
$y_i \in \Bbb A^d$ must coincide and the $t_i$'s must be equal to each
other. So our orbit contains the point $\frac{1}{q}gg_1x_1 + \dots +
\frac{1}{q}gg_qx_q$ of $K\# \dots \# K$.

\smallpagebreak

CASE $ r(P) > q$.  If the projections of the $|P|$ vectors $\{
g_\alpha (P) \bullet s_\alpha : \alpha \in [P] \}$ on the
$N$-dimensional vector space $\Bbb L^\perp$ are in general position,
$\sum_{\alpha \in [P]} t_\alpha (P)(g_\alpha (P) \bullet s_\alpha
)^\perp = 0$ is possible only when $|P| = N+1$, which we'll assume
from here on.

We will denote the $k$th coordinate of $s_\alpha \in S \subset \Bbb
A^d \subset \Bbb R^{d+1}$ by $s_{k,\alpha}$, and by $D(P,S)$ or just
$D$ the $N \times N$ determinant $|{g_\alpha}^ls_{k,\alpha}|$, with $1
\leq \alpha \leq N$ indexing the columns, and the rows indexed, in
lexicographic order, by the ordered pairs $(k,l)$, $1 \leq k \leq
d+1$, $1 \leq l \leq q-1$. Grouping together terms involving the first
$q-1$ rows, then the next $q-1$ rows, etc., we see that it has Laplace
expansion

$$ D = \sum_\pi (-1)^\pi \{ \Delta (g_{\pi_1}, \dots
,g_{\pi_{q-1}})s_{1,\pi_1} \cdots s_{1,\pi_{q-1}}\}\cdot \{ \cdots $$
$$ \cdots \}\cdot \{\Delta (g_{\pi_{N-q+1}}, \dots ,g_{\pi_N})
s_{d+1,\pi_{N-q+1}} \cdots s_{d+1,\pi_N}\},$$

\flushpar where $\pi$ runs over all permutations of $\{ 1, \dots , N \}$ such that $\pi_1 < \dots < \pi_{q-1}$; $\pi_q < \dots < \pi_{2(q-1)}$; etc.  So using the conjugation rule of Vandermonde determinants, we see that

$$\overline{D} = (-1)^{(d+1)\binom{q-1}{2}}(g_1 \cdots g_N)^{-q}D.$$

More generally we'll denote by $D_{m,N+1}$, $1 \leq m \leq N+1$, the
determinant obtained by replacing the $m$th column of $D$ by
$[{g_{N+1}}^ls_{k,N+1}]$ (so $D_{N+1,N+1}=D$) and use their Laplace
expansions and conjugation rules too.

We'll equip $\Bbb L$ with the \it Fourier basis \rm consisting of the
matrices $M_{l,k}$ : all rows zero except the $k$th which equals
$(\frac{1}{q}, \frac{1}{q}\omega ^l, \dots , \frac{1}{q}\omega
^{l(q-1)})$. Note that $(g_\alpha \bullet s_\alpha )^\perp = \sum \Sb
1\leq l < q \\ 1 \leq k \leq d+1 \endSb {g_\alpha}^ls_{k,\alpha}
M_{l,k}$; so our general position hypothesis is equivalent to saying
that \it all \rm the determinants $D_{m,N+1}$ are nonzero. So (this
step works even if \it any one \rm of the determinants is nonzero) the
real solution $<t_\alpha >$ of the $N$ equations $\sum_{\alpha
=1}^{N+1}t_\alpha {g_\alpha}^ls_{k,\alpha} = 0$ obeys $t_m = -
\frac{D_{m,N+1}}{D}t_{N+1}$ for $1 \leq m \leq N$. But these
determinantal ratios are real only if $(\frac{g_m}{g_{N+1}})^q = 1 \
\forall m$ which is impossible because there are $r(P) > q$ distinct
$g_\alpha$'s. \it q.e.d.

\medpagebreak

\rm REMARK 1. \it The above determinant $D(P,S)$ is zero for all $S \subset \Bbb A^d \subset \Bbb R^{d+1}$ if and only if $P \in S^{2N+1}$, $|P| = N+1$, is such that a $g_\alpha(P)$ repeats more than $d+1$ times for $1\leq \alpha \leq N$ \rm  (likewise for the other determinants $D_{m,N+1}$).  \lq If\rq \ is obvious, and for  \lq only if\rq \  note that if no $g_\alpha$ repeats more than $d+1$ times, then the $N = (d+1)(q-1)$ $\alpha$'s can be partitioned off into $d+1$ parts, of cardinality $q-1$ each, such that the $g_\alpha$'s corresponding to each part are distinct.  Assume without loss of generality that $\{ g_1, \dots ,g_{q-1} \}$ are all distinct, $\{ g_q , \dots g_{2(q-1)} \}$ are all distinct, and so on.  Now, choosing  $s_1 = \dots = s_{q-1} = (1,0,\dots ,0)$, $s_q = \dots = s_{2(q-1)} = (0,1,0, \dots ,0)$ and so on, we get  $D(P,S) = \Delta (g_1, \dots ,g_{q-1}) \cdots \Delta (g_{N-q+1}, \dots , g_N) \neq 0$. \it q.e.d. \rm

\medpagebreak

We'll denote by $\Bbb L_{\Bbb R}$ and $\Bbb L_{\Bbb R}^\perp$ the real
subspaces of $\Bbb L$ and $\Bbb L^\perp$ consisting of all \it real
matrices\rm . Note that a matrix $\lambda \in \Bbb L$ is real iff its
Fourier components $\lambda_{l,k}$, with respect to the basis $\{
M_{l,k} \}$, satisfy the conditions $\lambda _{l,k} =
\overline{\lambda_{q-l,k}}$. Indeed, as follows at once from $\lambda
= \sum_{l,k}\lambda_{l,k}M_{l,k}$ and $M_{l,k} =
\overline{M_{q-l,k}}$, the complex conjugate $\overline{\lambda}$ of
any $\lambda \in \Bbb L$ has Fourier components
$\overline{\lambda}_{l,k} = \overline{\lambda_{q-l,k}}$.  Note that
these real subspaces are preserved only by the $\Bbb Z/q$ action, and
likewise that complex conjugation commutes with the $\Bbb Z/q$, but
not with the $S^1$-action. Generalizing an argument used in the proof
of the previous lemma we'll now check the following.

\medpagebreak

\bf Lemma 3. \it A point $P\in S^{2N+1}$ lies on an $S^1$-orbit passing through the $q$-fold deleted join  if and only if there exist $v_\alpha \in \Bbb L^\perp _{\Bbb R}$, $\alpha \in [P]$, with \rm rank$\{ g_\alpha (P)\bullet v_\alpha  : \alpha \in [P]\} = |P| - 1$ \it and $0$ lies in the open convex hull of $\{ g_\alpha (P)\bullet v_\alpha : \alpha \in [P]\}$.

\medpagebreak

Proof. \rm Given any $P \in S^{2N+1}$ we can certainly find $|P|$
affinely independent vectors $w_\alpha \in \Bbb L^\perp _{\Bbb R}$
such that $0$ lies in their open convex hull.  (Moreover, these
vectors can be so chosen that $0$ has any prescribed positive
barycentric coordinates $t_\alpha$, say $t_\alpha = t_\alpha (P)$,
with respect to them.) When $P\in |K\ast \cdots \ast K|$, i.e. when
$g_\alpha (P) \in \Bbb Z/q$, we now get the required $v_\alpha \in
\Bbb L^\perp _{\Bbb R}$ by solving $w_\alpha = g_\alpha (P)\bullet
v_\alpha$. The same $v_\alpha$'s will work also for any other point
$gP$ of the orbit $\overline{P}$ through $P$.

Conversely, writing the given convex dependency $\sum_{\alpha \in
[P]}t_\alpha(g_\alpha (P)\bullet v_\alpha) = 0$ in components, we get
$\sum_{\alpha \in [P]} t_\alpha {g_\alpha}^lv_{\alpha,l,k} = 0$, $1
\leq k \leq d+1$, $1 \leq l < q$, i.e. $N$ equations in $|P|$ positive
variables $<t_\alpha >$. Since $v_\alpha \in \Bbb L^\perp _{\Bbb R}$,
multiplying the $\alpha$th column of the coefficient matrix
$[{g_\alpha}^lv_{\alpha , l,k}]$ by ${g_\alpha}^{-q}$ replaces the
element ${g_\alpha}^lv_{\alpha ,l,k}$ of its $(k,l)$th row by
${g_\alpha}^{l-q}v_{\alpha ,l,k}$ which is the complex conjugate of
the element ${g_\alpha}^{q-l}v_{\alpha ,q-l,k}$ of its $(k,q-l)$th
row. This gives the conjugation rule
$$ (\prod_{\alpha \neq \beta}(g_\alpha )^{-q}).D_{R,\beta} = \pm
\overline{D_{JR,\beta}}\ ,$$ where the sign is independent of $\beta$,
$D_{R,\beta}$ denotes determinant obtained by omitting the $\beta$th
column and using a cardinality $|P|-1$ subset $R$ of rows $(k,l)$, and
$JR$ denotes the corresponding subset of rows $(k,q-l)$. By hypothesis
we can choose $R$ so that at least one of these determinants
$D_{R,\beta}$ which we call $D_R$ is nonzero. Then the corresponding
$D_{JR}$ from the $D_{JR,\beta}$'s is also nonzero.  We now apply
Cramer's rule, to the rows $R$, and also to the rows $JR$, to solve
for $|P|-1$ of the $t_\alpha$'s in terms of one of them. This shows
that we must have $\frac{ D_{R,\beta}}{D_R} = \frac{
D_{JR,\beta}}{D_{JR}}$ and that these ratios are real. Applying the
above conjugation rule it now follows that all ratios of $ g_\alpha
(P)$'s must be in $\Bbb Z/q$, i.e. that $P$ must be on an $S^1$-orbit
passing through $|K\ast \cdots \ast K|$. \it q.e.d. \rm

\medpagebreak

Our map $s : S^{2N+1}\rightarrow \Bbb L^\perp$ was defined starting
from $\{ s_\alpha ^\perp \}$ which was the projection of a subset
constrained to be in $\Bbb A^d$. Taking our cue from the above lemma
it is better to consider \it the class of all subsets $\{ v_\alpha : 1
\leq \alpha \leq N+1\}$ of $\Bbb L^\perp _{\Bbb R}$. \rm This gives us
more \lq\lq elbow room\rq\rq \ e.g. by perturbing the set slightly we
can now assume (as we do in the following) that, for all $g_\alpha \in
\Bbb Z/q$, all minors of the $N \times (N+1)$ matrices $[g_\alpha
\bullet v_\alpha ]$ are nonzero.

\medpagebreak 

\bf Lemma 4. \it For any  subset $\{ v_\alpha \} \subset \Bbb L^\perp _{\Bbb R}$ which is generic in the above sense, there exists  an $S^1$-map $f : S^{2N+1}\rightarrow \Bbb L^\perp$, with  zeros only on orbits passing through the deleted join, and such that
$$ f(P) = \sum_{\alpha \in [P]} t_\alpha (P)( g_\alpha (P)\bullet
v_\alpha )$$ for all $P\in K\ast \cdots \ast K$.

\medpagebreak

Proof \rm \footnote{Sketchy ... needs clarifying and simplifying ...
comments most welcome ... }. We'll only use the weaker hypothesis that
rank$\{ g_\alpha \bullet v_\alpha : 1 \leq \alpha \leq N+1\} = N$ for
all $g \in \Bbb Z/q$, however some steps become slightly simpler under
the stronger assumption that all minors of these matrices are nonzero.

The basic idea in constructing $f$ is to use a \it subset $\{ f_\alpha
(P)\} \subset \Bbb L^\perp _{\Bbb R}$ \lq\lq moving\rq\rq\ with $P \in
S^{2N+1}$ \rm which coincides with the given \lq\lq constant\rq\rq\
$\{ v_\alpha \}$ on the deleted join. Here the $f_\alpha (P)$'s are
continuous and invariant under the $S^1$-action : $f_\alpha (gP) =
f_\alpha (P)\ \forall g \in S^1$. This ensures that the same recipe,
$$ f(P) = \sum_{\alpha \in [P]} t_\alpha (P)( g_\alpha (P)\bullet
f_\alpha (P)),$$ continues to define a continuous $S^1$-map
$S^{2N+1}\rightarrow \Bbb L^\perp$.

We'll equip the base space of our Hopf fibration $\pi
:S^{2N+1}\rightarrow \Bbb CP^N$, $P\mapsto \overline{P}$, with the \it
toric subdivision \rm $\Bbb CP^N = \cup_{\sigma \in K} \Bbb
CT^{\sigma}$, $\Bbb CT^{\sigma} = \{ \overline{P} : [P] = \sigma
\}$. Note that each $\Bbb CT^\sigma$ is homeomorphic to a product of
$|\sigma |-1$ copies of $\Bbb C^{\times}$, and that the fibration is
trivial over this \lq\lq complex torus\rq\rq . In fact for each
$\gamma \in \sigma$ the condition $g_\gamma (P)\equiv 1$ fixes a
continuous right inverse $\gamma :\Bbb CT^{\sigma}\rightarrow
S^{2N+1}$ of $\pi$. Let $\Bbb C^N_t = \cup_{|\sigma |\leq t}\Bbb
CT^\sigma$. The extension $f$ will be constructed inductively over
these skeletons using the fact that an $S^1$-map $S^{2N+1}\rightarrow
\Bbb L^\perp$ is same thing as a section of the bundle $\frak L^\perp
\rightarrow \Bbb CP^N$.

The main thing is to ensure at each step that the $f_\alpha$'s satisfy
the rank condition of Lemma 3 for then the reality of the $f_\alpha$'s
would imply that the $S^1$-map $f$ has no new zeros. This process is
easier to indicate for the \it case $q$ even \rm which we'll assume
now. Note that the $(k,l)$th component of $g_\alpha (P)^{-q/2}g_\alpha
(P)\bullet f_\alpha (P)$ is $g_\alpha^{-q/2}g_\alpha^l f_{\alpha ,l,k}
$ which is the complex conjugate of its $(k,q-l)$th component
$g_\alpha^{-q/2}g_\alpha^{q-l} f_{\alpha ,q-l,k} $. So like $f_\alpha
(P)$ this vector is also in $\Bbb L^\perp _{\Bbb R}$. Moreover clearly
any minor of $[g_\alpha (P)^{-q/2}g_\alpha (P)\bullet f_\alpha (P)]$
is nonzero iff the corresponding minor of $[g_\alpha (P)\bullet
f_\alpha (P)]$ is nonzero.

So we can work with the space of all $N\times (N+1)$ \lq\lq
real\rq\rq\ (but written in Fourier component) matrices $[ w_\alpha :
w_\alpha \in \Bbb L^\perp _{\Bbb R}]$ of rank $N$. Now this space has
the homotopy type of $GL_+ (N+1,\Bbb R)\simeq SO(N+1, \Bbb R)$ which
is fairly complicated. But we don't really have to deal with
obstruction-theoretic problems arising from this because to start with
we only have \it finitely many \rm such matrices $[ g_\alpha \bullet
v_\alpha : g_\alpha \in \Bbb Z/q, 1 \leq \alpha \leq N+1]$ and so we
can work within an open \it contractible \rm subspace $\Omega$ of rank
$N$ matrices of size $N\times (N+1)$ containing these finitely many
initial ones. (For Theorem 3 a slightly different Lemma 4 would be
needed because the set of initial matrices is not finite but they are
still defined over a space -- the union of the open top-most simplices
of the deleted join -- having the homotopy type of finitely many
points.) Restricting to columns indexed by $\alpha \in \sigma$ we get
a space of size $N \times |\sigma |$ matrices which we'll denote by
$\Omega _\sigma$. Note that these have rank at least $|\sigma
|-1$. We'll denote by $\Omega _\sigma ^0$ the subspace of those having
rank $|\sigma |$.

We'll assume that our inductive construction has taken care that over
any $\Bbb CT^\sigma$ with $|\sigma |\leq t$ the matrices $[g_\alpha
(\gamma \overline{P})^{-q/2}g_\alpha (\gamma \overline{P})\bullet
f_\alpha (\overline{P})]$ are in $\Omega _\sigma$ (or in
$\Omega^0_\sigma$ if we start with the all-minors-nonzero
hypothesis). Over the boundary of a $\Bbb CT^\theta$ with $|\theta | =
t+1$ we perturb the $f_\alpha (P)$'s to raise, if possible, the rank
of the similar matrix from $|\theta |-2$ to $|\theta |-1$ (or from
$|\theta |-1$ to $|\theta |$ if doing construction with the stronger
hypothesis on the $v_\alpha$'s). This is possible as long as there is
a row available outside a biggest sized nonzero minor because then by
perturbing the $f_{\alpha ,l,k}$'s of this row we can raise the rank.

Using the contractibility of $\Omega _{\theta}$ we now extend this
matrix $[g_\alpha (\gamma \overline{P})^{-q/2}g_\alpha (\gamma
\overline{P})\bullet f_\alpha (\overline{P})]$ to a matrix $A_\alpha
(\overline {P}) \in \Omega _{\theta}$ defined continuously over all of
$\Bbb CT^\theta$. Solving $A_\alpha (\overline{P})= g_\alpha (\gamma
\overline{P})^{-q/2}g_\alpha (\gamma \overline{P})\bullet f_\alpha
(\overline{P})$ we then get the required extensions of the functions
$f_\alpha$ over $\Bbb CT^\theta$. Using these extend the $S^1$-map
$f$. Applying Lemma 3 (with the stronger hypothesis one needs to use
this lemma only when $|\theta | = N+1$) it follows that no new zeros
have been introduced in this process.  \it q.e.d. \rm

\medpagebreak

REMARK 2. We note that for $q$ even the action $ [w_\alpha ]\mapsto
[g^{-q/2}g\bullet w_\alpha ]$ of $ g \in S^1$ on matrices multiplies
each minor by a nonzero constant and preserves the value of all
$N\times N$ minors because $(g^{-q/2})^N(g^1g^2...g^{q-1})^{d+1} =
1$. We note also that for the above argument it was necessary to use
\lq\lq moving\rq\rq\ points because, for $q \geq 3$, \it an $f$
defined by a \lq\lq constant\rq\rq\ $\{ v_\alpha \} \subset \Bbb
L^\perp _{\Bbb R}$ can not satisfy the required rank condition at all
$P$. \rm To see this, take any $1 \leq k \leq d+1$, and inductively
choose $g_\alpha$ as follows : in case $v_{\alpha ,1,k} =
\overline{v_{\alpha , q-1, k}} = 0$ or if there is no $\beta < \alpha$
with $v_{\beta ,1,k} = \overline{v_{\beta , q-1, k}}$ nonzero take any
$g_\alpha$, otherwise pick such a $\beta < \alpha$ and solve $\vmatrix
{g_\beta}v_{\beta ,1,k} & {g_{\alpha }}v_{\alpha
,1,k}\\{g_\beta}^{q-1}v_{\beta ,q-1,k} & {g_{\alpha }}^{q-1}v_{\alpha
,q-1,k}\endvmatrix = 0$ to find $g_\alpha \in S^1$. This implies that
all the $2\times 2$ determinants formed from the $(k,1)$th and the
$(k,q-1)$th rows of the $N\times N$ determinants $
|{g_\alpha}^lv_{\alpha ,l,k}|_{\alpha \neq m}$ are zero. So, at any
$P$ having such $g_\alpha$'s, these $N + 1$ determinants are all
zero. \it q.e.d. \rm

\heading 4. Proof of Theorem 1 \endheading 

If $S \subset \Bbb A^d$ has $t$ Tverberg partitions then any
sufficiently close $p(S) \subset \Bbb A^d$ has $\leq t$ Tverberg
partitions.  This follows because if the disjoint partition $\sigma_1
\cup \dots \cup \sigma_q$ of $S$ is such that $\text{conv}(\sigma_1)
\cap \dots \cap \text{conv}(\sigma_q) = \emptyset$ then the same is
true for the corresponding partition $p(\sigma_1) \cup \dots \cup
p(\sigma_q)$ of $p(S)$.  So without loss of generality we can assume
that $S \subset \Bbb A^d \subset \Bbb R^{d+1}$ is \it generic\rm ,
i.e. that the coordinates $s_{k,\alpha}$, $1 \leq k \leq d+1$, of its
points $s_\alpha$ are algebraically independent, but for the obvious
relations $\sum_{k=1}^{d+1}s_{k,\alpha} = 1$, over the field of
rationals.

Now consider the restriction $s^q : K \# \cdots \# K \rightarrow \Bbb
A^d \times \cdots \Bbb A^d$ of the $q$-fold cartesian product of the
linear map $s : |K| \rightarrow \Bbb A^d$ determined by $S$.  We note
that $K \# \cdots \# K$ is $(q-1)d$-dimensional, and can be equipped
with the cell subdivision provided by all cells of the type
$|\sigma_1| \times \cdots \times |\sigma_q|$, where the $\sigma_i$'s
are any pairwise disjoint nonempty subsets of $\{ e_1, \dots
,e_{N+1}\}$. Since $S$ is in general position on $\Bbb R^{d+1}$, the
$s^q$-images of the lower dimensional cells miss the diagonal of $\Bbb
A^d \times \cdots \times \Bbb A^d$, while the images of the top
dimensional cells either miss the diagonal, or else cut it cleanly in
just one interior point.

This last happens if and only if $(\sigma_1, \dots ,\sigma_q)$ is one
of the $q!$ permutations of a Tverberg partition $\{ \sigma_1, \dots
,\sigma_q \}$ of $S$.  Thus the number of points of $K \# \cdots \# K$
which are mapped to the diagonal by $s^q$ is $q!$ times the number of
Tverberg partitions of $S$, and so the number of $\Bbb Z/q$ orbits of
the $\Bbb Z/q$-space $K \# \cdots \# K$ which are imaged to a single
point is $(q-1)!$ times the number of Tverberg partitions of $S$.

We extend the $\Bbb Z/q$-map $s^q : K \# \cdots \# K \rightarrow \Bbb
A^d \times \cdots \times \Bbb A^d$ to the $S^1$-map $s$ of \S 3 and
consider its direct summand $s : S^{2N+1} \rightarrow \Bbb L^\perp$.
Identifying this $S^1$-map with the section $\Bbb CP^N @>s>> \frak
L^\perp$, $\overline{P} \mapsto [P,s(P)] \in (S^{2N+1}\times \Bbb
L^\perp)/S^1 = \frak L^\perp$, of the vector bundle $\frak L^\perp
\rightarrow \Bbb CP^N$, we see by using Lemma 2 that on $Q =
\{\overline{P} : r(P) \leq q \}$ this section has only \it Tverberg
zeros \rm $x$, i.e. $x = \overline{P}$ for some $P \in K\# \cdots \#
K$ which is mapped to the diagonal by $s^q$. The following lemma shows
that these zeros are isolated, and so that the section $s$ has no
other zeros in a sufficiently small neighbourhood of $Q$.

\medpagebreak

\bf Lemma 5. \it For  $S \subset \Bbb A^d$ generic, the section $s$  of $\frak L^\perp$ intersects its zero section $z$ transversely at all Tverberg zeros $x$. \rm  (Here as usual \lq\lq transversely\rq\rq \ means that the tangent space of $\frak L^\perp$ is the sum of the subspaces tangent to the two sections at the intersection.) \it 

\medpagebreak

Proof.  \rm Let $x = \overline{P}$ where $P \in K\# \cdots \# K$ is
mapped to the diagonal by $s^q$.  We must have $|P| = N+1$, for
otherwise, a proper subset of $S$ would have a Tverberg partition.
For the same reason (use Carath\'eodory's theorem) all Tverberg parts
$\sigma_i$ must have cardinality $\leq d+1$. Hence, by Remark 1, the
$N\times N$ determinants $D_{m,N+1}(P,S) =
|{g_\alpha}^ls_{k,\alpha}|_{\alpha \neq m}$ must be nonzero
polynomials in the $s_{k,\alpha}$'s.

The definitions $\overline{P} \mapsto [P,s(P)]$ and $\overline{P}
\mapsto [P,0]$ of $s$ and $z$ show that our assertion is equivalent to
checking that the map $\overline{P} \mapsto s(P)$, $g_{N+1}(P)=1$, has
a nonsingular jacobian at $x$. To verify this we'll use, as local
coordinates near $x$ the $2N$ reals $<t_\alpha , \theta_\alpha >$, $1
\leq \alpha \leq N$, where $g_\alpha = \text{cos}\theta _\alpha +
i\text{sin}\theta_\alpha$ .  In $\Bbb L^\perp$ we'll use as local
coordinates the $2N$ reals given by the real and imaginary parts of
the $(l,k)$th components, $1 \leq l < q$, $1 \leq k \leq d+1$, with
respect to the Fourier basis.  In these local coordinates our map
reads

$$ <t_\alpha , \theta_\alpha > \mapsto <s_{k,N+1}+\sum_\alpha t_\alpha
\{ \text{cos}(l\theta_\alpha )s_{k,\alpha} - s_{k,N+1}\} , \sum_\alpha
t_\alpha \text{sin} (l\theta_\alpha)s_{k,\alpha}>.$$

Computing partial derivatives with respect to the $t_\alpha$'s and
$\theta_\alpha$'s we obtain the jacobian $2N \times 2N$ determinant

$$
 \vmatrix
\text{cos}(l\theta_\alpha )s_{k,\alpha }- s_{k,N+1} & -lt_\alpha \text{sin}(l\theta_\alpha )s_{k,\alpha}\\
\text{sin}(l\theta_\alpha )s_{k,\alpha} & lt_\alpha \text{cos} (l\theta_{\alpha})s_{k,\alpha}
\endvmatrix  .
$$ Pulling the nonzero $t_\alpha$'s out from the columns, and doing
some row transformations we see that it is
$(\frac{1}{2})^N\prod_{\alpha = 1}^N t_\alpha$ times

$$ \vmatrix {g_\alpha}^ls_{k,\alpha} - s_{k,N+1} &
l{g_\alpha}^ls_{k,\alpha}\\ -\overline{g_\alpha}^ls_{k,\alpha}+
s_{k,N+1} & l\overline{g_\alpha}^ls_{k,\alpha}\\
\endvmatrix .
$$
Since, at $x$,\ $(g_\alpha)^q = 1 \ \forall \alpha$, on adding the
$(k,l)$th bottom row of this determinant to its $(k,q-l)$th top row,
we can change this to

$$ \vmatrix {g_\alpha}^ls_{k,\alpha}- s_{k,N+1} &
l{g_\alpha}^ls_{k,\alpha}\\ 0 & q\overline{g_\alpha}^ls_{k,\alpha}\\
\endvmatrix ,
$$
which equals $(D-D_{1,N+1}- \dots - D_{N,N+1})\cdot q^N\cdot
\overline{D} \neq 0$, because the left side is a nonzero polynomial in
the generic $s_{k,\alpha}$'s with coefficients in the algebraic number
field $\Bbb Q[\omega ]$.  \it q.e.d. \rm

\medpagebreak

Continuing the proof of Theorem 1 we now equip $\Bbb CP^N$, as well as
the complex $N$-dimensional fibres of $\frak L^\perp$, with the
orientations determined by their complex structures.  With respect to
these orientations one can speak of the \it local degree \rm $d_x$ of
the section $s$ of $\frak L^\perp \rightarrow \Bbb CP^N$ at each of
its isolated zeros $x$, i.e. the degree of the obvious map from the
link of $x$ to the unit sphere of the fibre $\frak L^\perp _{x}$ which
is determined by $s$.  From the above lemma it follows that this map
is a diffeomorphism, so this degree is $+1$ or $ -1$, depending on
whether or not the diffeomorphism preserves or reverses orientation.

We now perturb $\{ s_\alpha^\perp \} \subset \Bbb L^\perp _{\Bbb R }$
to a neighbouring $\{ v_\alpha \} \subset \Bbb L^\perp _{\Bbb R }$ for
which all minors of the matrices $[g_\alpha \bullet v_\alpha :
g_\alpha \in \Bbb Z/q ]$ are nonzero, and replace $s$ by the section
$f$ supplied by Lemma 4. This only perturbs the existing zeros
slightly (and they'll still be orbits passing through the deleted join
but maybe not the the deleted product) and introduces no new ones. At
each of the perturbed zeros one has the same local degree $\pm 1$ as
at the corresponding original one.

We recall now (this follows from the obstruction theoretic definition
of the Euler class) the well-known \it Poincar\' e-Hopf theorem \rm :
the sum $\sum_x d_x$ of the local degrees of $f$ coincides with the
Euler number $<e(\frak L^\perp ), \Bbb CP^N >$.  Here $e(\frak L^\perp
) = c_N(\frak L^\perp ) = c_N(\frak L )$, so this number is
$((q-1)!)^{d+1}$ by Lemma 1.  Since each $d_x = \pm 1$ it follows that
the number of Tverberg zeros is at least $((q-1)!)^{d+1}$, and thus
the number of Tverberg partitions of $S$ is at least $((q-1)!)^d$. \bf
q.e.d. \rm

\medpagebreak

REMARK 3.  An interesting \lq\lq stability\rq\rq \ deserves to be
mentioned : \it in all of the above we could have assumed $d > > q$.
\rm To see this, think of $\Bbb A^d$ as an affine hyperplane of $\Bbb
A^{d+1}$, perturb one of the points $v$ of the given general position
$S \subset \Bbb A^d$ into $\Bbb A^{d+1}$, and take $q-1$ new points on
the same side of $\Bbb A^{d+1}$.  For a Tverberg partition of this
subset of $\Bbb A^{d+1}$, the part which contains $v$ cannot contain a
new point; for then, one of the other parts will contain no new point,
implying that the common point of the convex hulls is on $\Bbb A^d$,
which is impossible for the proper subset $S \backslash \{ v\} $ has
(due to general position) no Tverberg partition.  So each Tverberg
partition of $S \subset \Bbb A^d$ can be associated with at most
$(q-1)!$ Tverberg partitions of the bigger set of $\Bbb A^{d+1}$, viz.
those obtained by adding one each of the new points to the $q-1$ parts
not containing $v$. Furthermore, if the new $q-1$ points are
coincident (or almost so) each of these is indeed a Tverberg partition
of this subset of $\Bbb A^{d+1}$.  \it q.e.d. \rm

\medpagebreak

It is also worth mentioning that the ratios of the coefficients
occuring in the expansion (cf. proof of Lemma 2) of $D(P,S)$, i.e.
the \it generalized cross ratios\rm ,
$$\frac {\Delta (g_{\pi_1},\dots ,g_{\pi_{q-1}})\cdots \Delta
(g_{\pi_{N-q+1}},\dots ,g_{\pi_N})}{\Delta (g_1,\dots ,g_{q-1})\cdots
\Delta (g_{N-q+1},\dots ,g_N)},$$ are always real. In fact they remain
the same under the stereographic projection $S^1 \backslash \{i\}
@>\cong>> \Bbb R$, $g \mapsto \text{cot}\frac{\theta}{2}$, shown
below.

$$\epsffile{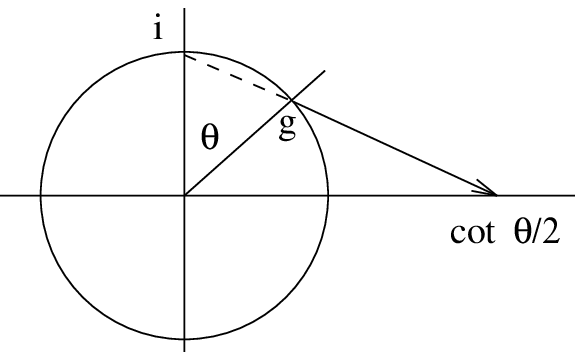}$$

\flushpar  To see this,  note that     $\overline{g} = g^{-1}$ implies that these ratios are indeed real, and then use  $|g_\alpha - g_\beta |\ = \ 2|\text{sin}\frac{\theta_\alpha}{2}\text{cos}\frac{\theta_\beta}{2} - \text{cos}\frac{\theta _\alpha}{2} \text{sin}\frac{\theta_\beta}{2} |$. The stereographic  projection also  shows that, when the  $g_\alpha (P)$'s are in  $\Bbb Z/q$, then the     $D_{m,N+1}(P,S)$'s are sort of like the cotangent or \lq\lq Dedekind sums\rq\rq \  of the field $\Bbb Q[\omega ]$, excepting that our  $s_{k,\alpha}$'s are not rational.

\medpagebreak   

Any cardinality $N+1$ set $S = \{ s_\alpha \} \subset \Bbb A^d$
assigns one such set of $N+1$ determinants to each top dimensional
simplex $(\sigma_1,\dots ,\sigma_q )$ of $K\ast \cdots \ast K$ :
simply use any point $P$ of this simplex. When $(d+1)\binom {q-1}{2}$
is even, these determinants are all real, otherwise all are purely
imaginary. We will also use the bigger determinant$\vmatrix 1\dots 1\\
{g_\alpha}^ls_{k,\alpha} \endvmatrix$, where $1\leq \alpha \leq N+1$
indexes the columns, and the rows, from the second onwards, are
indexed in lexicographic order by the pairs $(k,l)$, $1 \leq k \leq
d+1$, $1 \leq l < q$. Using the set $S \subset \Bbb A^d$ we now define
an $N$-dimensional \it characteristic cocycle \rm ${\Cal
X}_S(\sigma_1,\dots ,\sigma_q ) \in \{ -1,0,+1\}$ of $K\ast \cdots
\ast K$ : we assign to each top simplex the sign of the real number
$(-1)^N\vmatrix 1\dots 1\\ {g_\alpha}^ls_{k,\alpha} \endvmatrix \cdot
\overline{|{g_\alpha}^ls_{k,\alpha}|}_{\alpha \leq N}$.

\medpagebreak

\bf Theorem 2. \it For any generic cardinality $N+1$ subset $S \subset \Bbb A^d$  one has
$$\sum {\Cal X}_S(\sigma_1,\dots ,\sigma_q) \ = \ ((q-1)!)^d,$$ where
the summation is over all Tverberg partitions $ \{ s( \sigma_1),\dots
,s(\sigma _q ) \}$ of $S$.

\medpagebreak

Proof. \rm Our proof of Theorem 1 had given us the formula $\sum_x d_x
= ((q-1)!)^{d+1}$ where $x$ runs over all Tverberg zeros of $s$. From
the definition of $d_x$ given there it is clear that it is $+1$ or $
-1$ depending on whether or not the jacobian of that locally defined
map $\overline{P} \mapsto s(P)$, with respect to the coordinates
mentioned, is positive or not. We recall also that this jacobian had
turned out to be
$$
(\frac{q}{2})^N \prod_\alpha t_\alpha (D - D_{1,N+1} - \cdots -
D_{N,N+1}) \overline{D} = (- \frac{q}{2})^N \prod_\alpha t_\alpha
\vmatrix 1\dots 1\\ {g_\alpha}^ls_{k,\alpha} \endvmatrix
\overline{D}.$$ So it follows that $d_x = {\Cal X}_S(\sigma_1, \dots
,\sigma_q)$ where $(\sigma_1, \dots ,\sigma_q)$ denotes the top
simplex of the deleted join which contains $P$. We don't need to take
$g_{N+1}$ here because if each $g_\alpha$ is multiplied by $\omega$
the $(k,l)$th row of the bigger determinant is multiplied by $\omega
^l$ and that of the smaller by $\omega ^{-l}$. Also zeros $x$
interrelated by a reordering of $S$ are related by composing s with an
orientation preserving self-diffeomorphism of $\Bbb CP^N$ and so have
same $d_x$. Dividing both sides of the formula by $(q-1)!$ we obtain
the stated result.  \it q.e.d. \rm

\medpagebreak

REMARK 4.  From the viewpoint of number theory it is the opposite case
of \it rational \rm $S \subset \Bbb L^\perp _{\Bbb R}$ which is most
interesting. As the proof shows an index formula still holds provided
the section $s$ cuts the zero section $z$ transversely and this
happens not only for a general position $S$ but in some other cases
also.  For example \it $\Cal X_{S_o} = 1$ on all Tverberg partitions
of Sierksma's configuration $S_0$, \rm i.e. when $s_1 = \cdots =
s_{q-1} = (1,0,\dots ,0)$; $s_q = \cdots = s_{2(q-1)} = (0,1,\dots
,0)$; and so on; with the very last point being $({\frac{1}{d+1}},
\dots ,{\frac{1}{d+1}})$.  For this note that at any Tverberg
partition, all but the first of the entries of the last column of
$\vmatrix 1\dots 1\\ {g_\alpha}^ls_{k,\alpha} \endvmatrix$ can be made
zero by multiplying this column by $d+1$ and adding all the other
columns to it.  So this bigger determinant is $(-1)^N\frac
{N+d+1}{d+1}$ times the smaller and nonzero determinant
$|{g_\alpha}^ls_{k,\alpha}|_{\alpha \leq N}$. \it q.e.d. \rm

 The index formula probably gives non-trivial identities between the
 Dedekind sums of $\Bbb Q[\omega ]$ for suitably chosen rational sets
 $S \subset \Bbb L_{\Bbb R}^\perp$ having more than $((q-1)!)^{d+1}$
 zero orbits.

\heading 5. Concluding remarks \endheading 

We'll give some more applications of the above ideas, followed by a
few comments on their evolution.

\medpagebreak

\bf (a)  \rm  CONTINUOUS MAPS.\  Let $s:K \rightarrow \Bbb A^d$ be any continuous map,  and  let us say that a $q$-tuple $(x_1, \dots ,x_q)$ of points of $|K|$  is a \it separated \rm $q$-tuple point of $s$ if $s(x_1) = \dots = s(x_q)$ and one can find pairwise disjoint subsets $\sigma_1 , \dots , \sigma_q$ of $\{ e_1, \dots , e_{N+1}\}$ such that $x_1 \in |\sigma_1|, \dots , x_q \in |\sigma_q|$.  Then one has the following \lq\lq continuous\rq\rq \ generalization of Tverberg's theorem.  For \it primes $q$ \rm this was established by B\'ar\'any-Shlosman-Sz\"ucs [1], 1981, with a simpler proof (of which the following is an \lq\lq $S^1$-version\rq\rq ) given later in Sarkaria [6].

\medpagebreak

\bf Theorem 3. \it A continuous map from a $(q-1)(d+1)$-dimensional simplex to $\Bbb A^d$ has a separated $q$-tuple point. 

\medpagebreak

Proof. \rm Consider the $q$-fold join

$$ s^{(q)} : K\ast \cdots \ast K \rightarrow \Bbb R^{d+1} \times
 \cdots \times \Bbb R^{d+1} = \Bbb L_{\Bbb R},$$ of the given
 continuous $K @>s>> \Bbb A^d$, i.e. the continuous $\Bbb Z/q$-map
 which images $P = \sum_{i=1}^{r(P)} t_i(P)(g_i(P)\bullet x_i) \in K
 \ast \cdots \ast K$ to

$$s^{(q)}(P) = \sum_{i=1}^{r(P)} t_i(P)(g_i(P)\bullet s(x_i)) =
\sum_{\alpha \in [P]} t_\alpha(P)(g_\alpha(P)\bullet s_\alpha (P)),$$
where $s_\alpha (P) \in \Bbb A^d$ is defined, whenever $\alpha \in
[P]$, by
$$
s_\alpha (P) = s(\frac{\sum_{g_\beta (P) = g_\alpha (P)} t_\beta (P)
e_\beta }{\sum_{g_\beta (P) = g_\alpha (P)} t_\beta (P)}).
$$
Replacing each $s_\alpha (P)$ by its component $(s_\alpha (P))^\perp$
in the above formulae (use $\Bbb L_{\Bbb R} = \Bbb L_{\Bbb R}^0 \oplus
\Bbb L_{\Bbb R}^\perp$)) gives the direct summand
$$s^{(q)} : K \ast \cdots \ast K \rightarrow \Bbb L^\perp_{\Bbb R}.$$

Assume, if possible, that $s$ has no separated $q$-tuple points,
i.e. that the last map has no zeros. The same will be true if we
perturb $\{ s_\alpha (P)^\perp \}$ to a neighbouring $\{ v_\alpha
(P)\} \subset \Bbb L^\perp_{\Bbb R}$ and use this to define our map
$K\ast \cdots \ast K \rightarrow \Bbb L^\perp _{\Bbb R}$.

Choosing these $v_\alpha (P)$'s so that with $g_\alpha \in \Bbb Z/q$
the matrices $[g_\alpha \bullet v_\alpha (P)]$ have minors nonzero we
can now extend this $\Bbb Z/q$-map, again using a Lemma 3 dependent
construction analogous to that of Lemma 4, to an $S^1$-map
$S^{2N+1}\rightarrow \Bbb L^\perp$, having no zeros anywhere on
$S^{2N+1}$. This contradicts the fact that the Euler class of $\frak
L^\perp$ being nonzero, it admits no everywhere nonzero continuous
section. \it q.e.d. \rm

\medpagebreak

 The Borsuk-Ulam theorem says that there is no continuous $\Bbb
 Z/2$-map from a free $\Bbb Z/2$-sphere to a lesser dimensional free
 $\Bbb Z/2$-sphere, and only essential ones to an equi-dimensional
 sphere.  Indeed it is known that the same is true for all finite
 groups $G \neq 1$. The $N$-dimensional free $\Bbb Z/q$-complex $K\ast
 \cdots \ast K$ being $(N-1)$-connected, this immediately gave [6] the
 required contradiction, \it for the case $q$ prime, \rm because then
 the $\Bbb Z/q$-action on $\Bbb L_{\Bbb R}^\perp \backslash \{ 0\}
 \simeq S^{N-1}$ is free . A similar argument also gave a sharp
 generalization of the Van Kampen-Flores theorem to all primes $q \geq
 2$. This approach also gives (for $q$ prime) an estimate -- see Vu\u
 ci\'c-\u Zivaljevi\' c [15] -- for the least number of Tverberg
 partitions which however falls short of the Sierksma number.

\medpagebreak

\bf (b) \rm PIECEWISE LINEAR MAPS. \  The main theorem generalizes to any subclass of continuous functions $s: K \rightarrow \Bbb A^d$ defined by a \lq\lq local smoothness or finiteness\rq\rq \ condition which is strong enough to guarantee transversality, and thus $d_x = \pm 1$, at the generic zeros, e.g. one has the following.

\medpagebreak

\bf Theorem 4. \it Any piecewise linear continuous map from a $(q-1)(d+1)$ dimensional simplex to $\Bbb A^d$ has at least $((q-1)!)^d$ separated $q$-tuple points.
\rm (Here as usual \lq\lq piecewise linear\rq\rq\ means that the map is linear with respect to some finite subdivision.)

\medpagebreak

\it Proof. \rm Let $s: K \rightarrow \Bbb A^d$ be linear on the subdivision $K'$ of $K$. Clearly the number of separated $q$-tuple points can not increase if $s$ is replaced by any neighbouring map also linear on $K'$. So we can assume that the set of coordinates of the values which $s$ takes on the vertices of $K'$ is generic in the sense of \S 4.

Using such an $s$ now define $S^1$-map $f : S^{2N+1} \rightarrow \Bbb
L^\perp$ as in the proof of Theorem 3. Considered as a section of
$\frak L^\perp$, the number of zeros of $f$ on the projection (in
$\Bbb CP^N$) of the $q$th deleted join will be $(q-1)!$ times the
number of separated $q$-tuple points. One can check, by a calculation
like that in the proof of Lemma 5, that at these points the section
$f$ cuts the zero section of $\frak L^\perp$ transversely.  Then the
same argument as in \S 4 applies and gives the required result. \it
q.e.d. \rm

\medpagebreak

\bf (c) \rm A GLOBALIZATION.\   We now let K denote any $N$-dimensional simplicial complex on the set of vertices $\{e_1, \dots ,e_M\}$, where $M \geq N+1$. Let $\text{Sph}(K) \subset S^{2M-1} \subset \Bbb C^M$ be the $S^1$-subspace consisting of all points $\sum_{\alpha \in [P]}t_\alpha (P)g_\alpha (P)e_\alpha$ with $[P] \subset \{e_1, \dots ,e_M\}$ a simplex of $K$. The space obtained by dividing this out by $S^1$ will be denoted $\text{Proj}(K)$, and an oriented  vector bundle $\frak L^\perp \rightarrow \text{ Proj}(K)$ defined on it exactly as before.

\medpagebreak

\bf Theorem 5. \it For any linear map $s$ from a $(q-1)(d+1)$ dimensional simplicial complex $|K|$ to $\Bbb A^d$ there exist at least \rm $<e(\frak L^\perp ), \text{ Proj}(K)> \div\  (q-1)!$ \it  partitions $\sigma_1 \cup \dots \cup \sigma_q$ of simplices of $K$ into pairwise disjoint faces such that \rm $\text{conv}(s\sigma_1)\cap \dots \cap \text{conv}(s\sigma_q)$ \it is nonempty.

\medpagebreak

\rm Here we note that $\text{Proj}(K)$ is a union of complex projective spaces and as such has a well defined fundamental class, which too has been denoted $\text{Proj}(K)$ in the above statement.  We omit the proof which is just as before; also it has a  piecewise linear generalization like Theorem 4, and can be sharpened to an index formula  like that of Theorem  2, with   characteristic cocycle of  a cardinality $M$ subset of $\Bbb A^d$ defined just as before. Note that  for    $q=2$  this  index formula says that the algebraic number of circuits  of the \lq\lq oriented matroid\rq\rq\  determined by this affine  set coincides  with the Euler number of a certain  vector bundle. 

\medpagebreak

  One should think -- for more on these lines see [9] and the lecture
  notes (under preparation) of my Panjab University Topology Seminar
  of 1994-95 -- of $\text{Sph}(K)$ as a \lq\lq visualisation\rq\rq ,
  just like the much better known $|K|$, of the simplicial complex
  $K$. Indeed there are many others, e.g. in a subsequent paper we'll
  give some interesting applications of $S^1$-versions of the \lq\lq
  bigger\rq\rq\ deleted joins used in Van Kampen's embedding
  theory. We note also that the incidence rule used in [9] to define a
  \lq\lq cyclotomic homology\rq\rq \ of simplicial complexes is like
  the definition of $\Cal X_S$ in \S 4.

\medpagebreak

\bf (d) \rm A LOOK BACK.\  Shortly after writing [6] I had made an attempt [7] to interpret the Sierksma number as a certain invariant sum of local degrees at the intersections of $s^q(K\# \cdots \# K)$ with the diagonal of $\Bbb A^d \times \cdots \times \Bbb A^d$. This was interesting -- e.g. a cycle over the cyclotomic field  $\Bbb Q[\omega ]$ was used -- but unsuccessful: a computation (some initially  overlooked  sign changes  were later pointed out to me by Kalai) using Sierksma's configuration $S_0$ shows that one only gets zero ! 

This had happened because a characteristic class of a \it finite \rm
group action was being evaluated on the aforesaid cycle; so it was
clear from that point on that an infinite group was needed, with $S^1$
being the obvious choice. (Intuitively this \lq\lq
complexification\rq\rq \ has the effect of changing the errant vertex
transpositions into vertex-pair transpositions, and thus a total
cancellation of the local degrees does not happen.)  The mechanics of
doing this became clear to me much later when I spoke about the
aforementioned \lq\lq visualisations\rq\rq \ in my 1994-95 topology
seminar, and the essentials of the above proof, including the key
Vandermonde conjugation argument of Lemma 2, were in hand by the end
of summer 1995. A more careful analysis -- comments of Ofer Gabber
were of great help in this regard -- however showed that we were still
short of a complete knowledge of all the zeros of the canonically
defined section $s$.  We have now side-stepped this difficulty by
using moving subsets $\{ f_\alpha (P)\} \subset \Bbb L^\perp _{\Bbb
R}$.

We note that Lemma 1 also follows, by virtue of the Poincar\'e-Hopf
theorem, from Remark 4 which showed that all the $((q-1)!)^{d+1}$
Tverberg zeros $x$ of Sierksma's configuration $S_0$ have the same
local degree $d_x = 1$. Thus our proof is in the same spirit as the
original inspiration of [7], viz. an argument given by Van Kampen in
his amazing paper [14] of 1932.  He showed that for $n\geq 2$ the
algebraical number of separated general position self-intersections,
of an $n$-complex immersed in $\Bbb A^{2n}$, is invariant under
deformations.  Then a computation, using a particular immersion of the
$n$-skeleton of a $(2n+2)$-simplex, shows that this complex can not
embed in $\Bbb A^{2n}$.

\medpagebreak

\bf (e) \rm WHAT LIES AHEAD ? \ We'll give in a sequel  more regarding the combinatorics  of the signs $d_x$,  some generalizations to skeletons of simplices, and some interesting applications of our index formula. Indeed, since other \it characteristic classes\rm , of say the pseudomanifolds $\text{Proj}(K)$, can be defined inductively in  terms of  suitable split  bundles, we hope that, like $e(\frak L^\perp )$, these too contain analogous combinatorial cocycles, and that there are  similar  combinatorial \it index formulas \rm for  other   numerical topological invariants.

We gave in [8] a very simple proof of Tverberg's theorem -- see
esp. Onn's remark (3), also see Kalai [4] -- which again uses, like
[6] and the proof of Theorem 2 above, the matrix space $\Bbb L_{\Bbb
R}^\perp$, but avoids topology by exploiting instead the linearity of
$s$ via an elementary convexity argument of B\'ar\'any. The above
proof of Theorem 1 can also probably be simplified to one using only
the \it representation theory \rm of $S^1$.  Also we're trying to make
another proof in which one replaces $S^1$ by $\Bbb C^\times$ and uses
\it field theory \rm : in this context we remark that apparently --
cf. Sullivan [12] -- the (discontinuous!) \lq\lq Galois
symmetries\rq\rq , of the $\Bbb Z/q^n$-deleted joins contained in
their covers, have much to say about the homotopy and homeomorphism
classification of the complex varieties $\text{Proj}(K)$. Finally, we
feel that \it noncommutative \rm versions of these arguments,
e.g. using $SU(2)$ instead of $S^1$, will be even more insightful.

\enddocument